\documentclass[draft,leqno]{amsart}

\usepackage[Symbol]{upgreek}

\usepackage{amssymb}
\usepackage{epic}
\usepackage{mathrsfs}
\usepackage{accents}
\usepackage{stmaryrd}

\input{xy}
\xyoption{matrix}
\xyoption{arrow}

\numberwithin{equation}{section}

\newcommand{\bbold}{\mathbb}

\def\R { {\bbold R} }
\def\Q { {\bbold Q} }
\def\Z { {\bbold Z} }
\def\Co { {\bbold C} }
\def\N { {\bbold N} }
\def\T { {\bbold T} }

\def \ex{\operatorname{e}}

\def \Frac {\operatorname{Frac}}

\renewcommand\epsilon{\varepsilon}

\def \<{\langle}
\def \>{\rangle}

\def \((  {(\!(}
\def \)) {)\!)}

\def \res{\operatorname{res}}

\DeclareMathSymbol{\precequ}{\mathrel}{symbols}{"16}
\DeclareMathSymbol{\succequ}{\mathrel}{symbols}{"17}

\def \nasymp{\not\asymp}

\newtheorem{theorem}{Theorem}[section]
\newtheorem{lemma}[theorem]{Lemma}
\newtheorem{prop}[theorem]{Proposition}
\newtheorem{cor}[theorem]{Corollary}

\newtheorem*{conjecture}{Conjecture}

\theoremstyle{definition}

\theoremstyle{remark}

\let\oldi\i
\let\oldj\j
\renewcommand\i{\relax\ifmmode{\boldsymbol{i}}\else\oldi\fi}
\renewcommand\j{\relax\ifmmode{\boldsymbol{j}}\else\oldj\fi}

\renewcommand\leq{\leqslant}
\renewcommand\geq{\geqslant}
\renewcommand\preceq{\preccurlyeq}
\renewcommand\succeq{\succcurlyeq}
\renewcommand\le{\leq}
\renewcommand\ge{\geq}
\renewcommand\frak{\mathfrak}

\DeclareMathAlphabet{\mathbf}{OML}{cmm}{b}{it}

\DeclareFontFamily{U}{fsy}{}
\DeclareFontShape{U}{fsy}{m}{n}{<->s*[.9]psyr}{}
\DeclareSymbolFont{der@m}{U}{fsy}{m}{n}
\DeclareMathSymbol{\der}{\mathord}{der@m}{182}

\DeclareSymbolFont{der@m}{U}{fsy}{m}{n}
\DeclareMathSymbol{\derdelta}{\mathord}{der@m}{100}

\DeclareSymbolFont{imag@m}{OT1}{cmr}{m}{ui}
\DeclareMathSymbol{\imag}{\mathord}{imag@m}{105}

\DeclareFontFamily{OMS}{smallo}{}
\DeclareFontShape{OMS}{smallo}{m}{n}{<->s*[.65]cmsy10}{}
\DeclareSymbolFont{smallo@m}{OMS}{smallo}{m}{n}
\DeclareMathSymbol{\smallo}{\mathord}{smallo@m}{79}

\DeclareFontFamily{OMS}{largerdot}{}
\DeclareFontShape{OMS}{largerdot}{m}{n}{<->s*[.8]cmsy10}{}
\DeclareSymbolFont{largerdot@m}{OMS}{largerdot}{m}{n}
\DeclareMathSymbol{\largerdot}{\mathord}{largerdot@m}{15}

\DeclareMathSymbol{\llambda}{\mathord}{der@m}{108}
\DeclareMathSymbol{\rrho}{\mathord}{der@m}{114}

\def \upg{\upgamma}
\def \Upg{\Upgamma}
\def \upl{\uplambda}
\def \Upl{\Uplambda}
\def \upo{\upomega}

\def \Upd{\Updelta}

\newcommand{\equationqed}[1]{\[\pushQED{\qed}#1 \qedhere\popQED\]\let\qed\relax}
\newcommand{\alignqed}[1]{\begin{align*}\pushQED{\qed} #1 \qedhere\popQED\end{align*}\let\qed\relax}

\makeatletter
\newcommand{\dminus}{\mathbin{\text{\@dminus}}}

\newcommand{\@dminus}{%
  \ooalign{\hidewidth\raise1ex\hbox{\bf.}\hidewidth\cr$\m@th-$\cr}%
}
\makeatother

\def \C{\mathcal{C}}
\def \Car{\mathcal{C}^r_a}
\def \Caz{\mathcal{C}^0_a}
\def \Cao{\mathcal{C}^1_a}
\def \Cat{\mathcal{C}^2_a}
\def \Cainf{\mathcal{C}^{\infty}_a}
\def \Caom{\mathcal{C}^{\omega}_a}

\def \Gr{\mathcal{C}^r}
\def \Gz{\mathcal{C}^0}
\def \Go{\mathcal{C}^1}
\def \Gt{\mathcal{C}^2}
\def \Gi{\mathcal{C}^{<\infty}}
\def \Ginf{\mathcal{C}^{\infty}}
\def \Gom{\mathcal{C}^{\omega}}
\def \inv{\operatorname{inv}}
\def \Sol{\operatorname{Sol}}

\begin{document}

\title{Hardy Fields, the Intermediate Value property, and $\upo$-Freeness}

\author[Aschenbrenner]{Matthias Aschenbrenner}
\address{Department of Mathematics\\
University of California, Los Angeles\\
Los Angeles, CA 90095\\
U.S.A.}
\email{matthias@math.ucla.edu}

\author[van den Dries]{Lou van den Dries}
\address{Department of Mathematics\\
University of Illinois at Urbana-Cham\-paign\\
Urbana, IL 61801\\
U.S.A.}
\email{vddries@math.uiuc.edu}

\author[van der Hoeven]{Joris van der Hoeven}
\address{\'Ecole Polytechnique\\
91128 Palaiseau Cedex\\
France}
\email{vdhoeven@lix.polytechnique.fr}

\begin{abstract} We discuss the conjecture that every maximal Hardy field has the Intermediate Value Property for differential polynomials, and its equivalence to the statement that all maximal Hardy field are elementarily equivalent to the differential field of transseries. As a modest but essential step towards establishing the conjecture we show that every maximal Hardy field is $\upo$-free.
\end{abstract}

\thanks{The first-named author was partially supported by NSF Grant DMS-1700439.}

\date{March 2019}

\maketitle

\section*{Introduction}\label{intro}

\noindent
Du Bois Reymond's ``orders of infinity'' were put on a firm
basis by Hardy~\cite{Ha} and Hausdorff~\cite{Hau}, leading to the notion of a Hardy field (Bourbaki~\cite{Bou}). A Hardy field is a field $H$ of germs at $+\infty$ of differentiable real-valued functions
on intervals~$(a,+\infty)$ such that the germ of the derivative
of any differentiable function whose germ is in $H$ is also in 
$H$. (See Section~\ref{bh} for more precision.) 
A Hardy field is naturally a differential field, and is
an ordered field with the germ of $f$ being $>0$ if and only if~$f(t)>0$, eventually. 

If $H$ is a Hardy field, then so is $H(\R)$ (obtained by adjoining the germs of the constant functions)
and for any $h\in H$, the germ $\ex^h$ generates a Hardy field~$H(\ex^h)$ over~$H$, and so does any differentiable germ with derivative
$h$~\cite{Ros}. Each Hardy field~$H$ has a unique Hardy field extension
that is algebraic over $H$ and real closed~\cite{Rob}. The ultimate extension result of this kind would be the following:

\begin{conjecture}
Let $H$ be a Hardy field, $P(Y)\in H\{Y\}$ a differential polynomial and~${f< g}$ in $H$ such that $P(f) <0 <  P(g)$. Then there is an element $\phi$ in a Hardy field extension
of $H$ such that $f < \phi < g$ and $P(\phi)=0$.
\end{conjecture}

\noindent
In \cite{D} this is proved for $P$ of order $1$.  In  \cite{vdH:hfsol} it is shown that there do exist Hardy fields with the intermediate value property for all differential polynomials. 
Every Hardy field extends to a maximal Hardy field, by Zorn, and
so the Conjecture above is equivalent to maximal Hardy fields having the intermediate value property for differential polynomials. By the results mentioned earlier, maximal Hardy fields contain $\R$  as a subfield and are Liouville closed in the 
sense of \cite{ADH}.  At the end of Section~\ref{hf} we show that for Liouville closed Hardy fields containing $\R$ the 
intermediate value property is equivalent to the conjunction of two other properties, {\em $\upo$-freeness\/} and {\em newtonianity}. These two notions are central in \cite{ADH} in a more general setting. Roughly speaking, $\upo$-freeness controls the solvability of second-order homogeneous linear differential equations in suitable extensions, and newtonianity is a very strong version of differential-henselianity. (We did not consider the intermediate value property in \cite{ADH} and mention it here
mainly for expository reasons: it  is easier to grasp than the more subtle and more fundamental notions of
$\upo$-freeness and newtonianity.)

The main result in these seminar notes is that any Hardy field has an $\upo$-free Hardy field extension: Theorem~\ref{upo}.  
We do not present it here, but we also have a detailed outline for showing  that any $\upo$-free Hardy field
extends to a newtonian $\upo$-free Hardy field. At this stage (March 2019) that proof is not yet finished. 
If we finish the proof, it would follow that every maximal Hardy field is an $\upo$-free newtonian Liouville closed $H$-field with small derivation, in the terminology of \cite{ADH}. Now 
the elementary theory $T^{\text{nl}}_{\text{small}}$ of $\upo$-free newtonian Liouville closed 
$H$-fields with small derivation is 
{\em complete}, by \cite[Corollary 16.6.3]{ADH}. Thus finishing the proof alluded to would give that any two maximal Hardy fields are indistinguishable as to their elementary properties, that is, any two maximal Hardy fields would be
elementarily equivalent as ordered differential fields.

\medskip\noindent
The present seminar notes prove some results announced in our exposition~ \cite{ADH2}. There we also
discuss another fundamental conjecture and partial results towards it, namely that the underlying ordered set of
any maximal Hardy field is $\eta_1$. Our plan for proving it  does depend on first establishing the conjecture that
maximal Hardy fields are newtonian. The two conjectures together imply:  all maximal Hardy fields are isomorphic as ordered differential fields, assuming the  continuum hypothesis (CH); for more on this, see~\cite{ADH2}.

\medskip\noindent
Let us add here a remark about maximal Hardy fields that is more set-theoretic in nature.  Every Hardy field is contained in
the ring $\C$ of germs at $+\infty$ of continuous real-valued functions on half-lines $(a,+\infty)$, and so there {\em at most\/} $2^{\frak{c}}$ many Hardy fields, where
$\frak{c}=2^{\aleph_0}$ is the cardinality of the continuum; note that $\frak{c}$ is  also the cardinality of $\C$. It is worth mentioning that the two conjectures above imply that there are in fact $2^{\frak{c}}$ many different
maximal Hardy fields: In an email to one of the authors, Ilijas Farah  showed that there are $2^{\frak{c}}$ many 
maximal Hausdorff fields, Hausdorff fields being the subfields of the ring $\C$ (without differentiability assumptions as in the case of Hardy fields). Farah's proof can easily be modified to give the same conclusion about the number of maximal Hardy fields assuming these conjectures.

\medskip\noindent
Throughout we use the
algebraic and valuation-theoretic tools from ~\cite{ADH}. We need in addition 
analytic facts about real and complex solutions of linear differential equations; these facts and various generalities about Hardy fields are in Section~\ref{bh}.   
 
\medskip\noindent 
The second-named author gave a talk about the above material
in the {\it S\'eminaire de structures alg\'ebriques ordonn\'ees}\/ in honor of Paulo Ribenboim's~90th birthday.
We dedicate this paper to Paulo in gratitude for his fundamental contributions to the theory of valuations, which is indispensable  in our work.

\subsection*{Notations and terminology} Throughout, $m$,~$n$ range over $\N=\{0,1,2,\dots\}$.
Given an additively written abelian group $A$ we let $A^{\ne}:=A\setminus\{0\}$.
 Rings are commutative with  identity $1$, and for a ring~$R$ we let $R^\times$ be the multiplicative group of units (consisting of the $a\in R$ such that $ab=1$ for some $b\in R$). A {\em differential ring\/} will be a ring $R$ containing
(an isomorphic copy of) $\Q$ as a subring and equipped with a derivation~$\der\colon R \to R$; note that then $C_R:=\big\{a\in R:\ \der(a)=0\big\}$ is a subring of $R$, called the ring of constants of $R$, and that $\Q\subseteq C_R$.  If $R$ is a field, then so is $C_R$.
An {\em ordered differential field}\/ is an ordered field  equipped with a derivation; such an ordered differential field is
in particular a differential ring.

Let $R$ be a differential ring and $a\in R$. When its derivation $\der$ is clear from the context we denote $\der(a),\der^2(a),\dots,\der^n(a),\dots$ by $a', a'',\dots, a^{(n)},\dots$, and if $a\in R^\times$, then $a^\dagger$ denotes $a'/a$, so $(ab)^\dagger=a^\dagger + b^\dagger$ for $a,b\in R^\times$. In Sections~\ref{hf} and~\ref{secupo} we need to consider the function $\omega=\omega_R\colon R \to R$ given by $\omega(z)=-2z'-z^2$, and
the function $\sigma=\sigma_R\colon R^\times \to R$ given by
$\sigma(y)=\omega(z)+y^2$ for $z:= -y^\dagger$.

We have the differential ring $R\{Y\}=R[Y, Y', Y'',\dots]$ of differential polynomials in an indeterminate $Y$. 
We say that  
$P=P(Y)\in R\{Y\}$ has order at most $r\in \N$ if
$P\in R[Y,Y',\dots, Y^{(r)}]$; in this case $P=\sum_{\i}P_{\i}Y^{\i}$, as in \cite[Section~4.2]{ADH}, with~$\i$ ranging over tuples
$(i_0,\dots,i_r)\in \N^{1+r}$, $Y^{\i}:= Y^{i_0}(Y')^{i_1}\cdots (Y^{(r)})^{i_r}$,     coefficients~$P_{\i}$   in $R$, and $P_{\i}\ne 0$ for only finitely many $\i$. For $P\in R\{Y\}$ and $a\in R$ we let $P_{\times a}(Y):=P(aY)$.
For $\phi\in R^\times$ we let $R^{\phi}$ be the {\it compositional conjugate of $R$ by~$\phi$}\/: the differential ring
with the same underlying ring as~$R$ but with derivation~$\phi^{-1}\der$ instead of $\der$. 
We have an $R$-algebra isomorphism $P\mapsto P^\phi\colon R\{Y\}\to R^\phi\{Y\}$ such that $P^\phi(y)=P(y)$ for all $y\in R$;
see \cite[Section~5.7]{ADH}.

\bigskip\noindent
For a field $K$ we have $K^\times=K^{\ne}$, and
a (Krull) valuation on $K$ is a surjective map 
$v\colon K^\times \to \Gamma$ onto an ordered abelian 
group $\Gamma$ (additively written) satisfying the usual laws, and extended to
$v\colon K \to \Gamma_{\infty}:=\Gamma\cup\{\infty\}$ by $v(0)=\infty$,
where the ordering on $\Gamma$ is extended to a total ordering
on $\Gamma_{\infty}$ by $\gamma<\infty$ for all 
$\gamma\in \Gamma$. 
A {\em valued field\/} $K$ is a field (also denoted by $K$) together with a valuation ring~$\mathcal O$ of that field,
and the corresponding valuation $v\colon K^\times \to \Gamma$  on the underlying field is
such that $\mathcal O=\{a\in K:va\geq 0\}$ as explained in~\cite[Section~3.1]{ADH}.

\medskip
\noindent
Let $K$ be a valued field
with valuation ring $\mathcal O_K$ and valuation $v\colon K^\times \to \Gamma_K$. Then~$\mathcal O_K$ is a local ring  with maximal ideal $\smallo_K=\{a\in K:va>0\}$  and residue field $\res(K)=\mathcal{O}_K/\smallo_K$. 
If $\res(K)$ has characteristic zero, then $K$ is said to be of equicharacteristic zero.
When the ambient valued field $K$ is clear from the context, then 
we denote $\Gamma_K$, $\mathcal O_K$, $\smallo_K$, by $\Gamma$, $\mathcal O$, $\smallo$, respectively,
and for $a,b\in K$ we set 
\begin{align*} a\asymp b &\ :\Leftrightarrow\ va =vb, & a\preceq b&\ :\Leftrightarrow\ va\ge vb, & a\prec b &\ :\Leftrightarrow\  va>vb,\\
a\succeq b &\ :\Leftrightarrow\ b \preceq a, &
a\succ b &\ :\Leftrightarrow\ b\prec a, & a\sim b &\ :\Leftrightarrow\ a-b\prec a.
\end{align*}
It is easy to check that if $a\sim b$, then $a, b\ne 0$, and that
$\sim$ is an equivalence relation on $K^\times$. 
We use {\em pc-sequence\/} to abbreviate
{\em pseudocauchy sequence}, and $a_\rho\leadsto a$ indicates that the pc-sequence $(a_\rho)$ pseudoconverges to $a$; see \cite[Sections~2.2,~3.2]{ADH}.
As in \cite{ADH},  a {\em valued differential field\/} is a valued field $K$ of equicharacteristic zero that is also equipped with a derivation 
$\der\colon K \to K$, and an 
 {\em ordered valued differential field\/} is a valued differential field $K$ equipped with an ordering on $K$ making $K$ an ordered field.

\section{$H$-Fields and IVP}\label{hf}

\noindent
We recall from \cite[Introduction]{ADH} that an {\it $H$-field}\/ is an ordered differential field $K$ with constant field~$C$ such that: \begin{enumerate}
\item[(H1)] $\der(a)>0$ for all $a\in K$ with $a>C$;
\item[(H2)] $\mathcal{O}=C+\smallo$, where $\mathcal{O}$ is the
convex hull of $C$ in the ordered field $K$, and $\smallo$ is the maximal ideal of the valuation ring $\mathcal{O}$.
\end{enumerate}
Let $K$ be an $H$-field, and let $\mathcal{O}$ and $\smallo$ be as in (H2). Thus $K$ is a valued field with valuation ring $\mathcal{O}$. 
The residue morphism $\mathcal O\to\res(K)=\mathcal O/\smallo$ restricts to an isomorphism $C\xrightarrow{\cong}\res(K)$.
The valuation topology on $K$ equals its order topology if $C\ne K$. We consider $K$ as an 
$\mathcal{L}$-structure, where 
$$\mathcal{L}\ :=\ \{\,0,\,1,\, {+},\, {-},\, {\times},\, \der,\, {<},\, {\preceq}\,\}$$ is the language of ordered valued differential fields. The symbols~$0$,~$1$,~${+}$,~${-}$,~${\times}$,~$\der$,~$<$ are interpreted as usual in $K$, and 
$\preceq$ encodes the valuation: for $a,b\in K$, 
$$ a\preceq b\quad \Longleftrightarrow\quad a\in \mathcal{O} b.$$ 
An $H$-field $K$ is said to be {\it Liouville closed}\/ if it is real closed and for all $a\in K$ there exists $b\in K$ with $a=b'$ and also $b\in K^\times$ with $a=b^\dagger$.

\subsection*{Remarks on IVP} Ordered valued differential subfields of $H$-fields are called pre-$H$-fields, and are characterized in \cite[Section~10.5]{ADH}.  Below we assume
some familiarity with the $H$-asymptotic couple $(\Gamma, \psi)$ of a pre-$H$-field $K$, as explained in \cite{ADH}, and
properties of $K$ based on those of $(\Gamma, \psi)$, such as
$K$ having {\it asymptotic integration}\/ and~$K$ having a {\it gap}\/ \cite[Sections~9.1, 9.2]{ADH}.

Let $K$ be a pre-$H$-field. We say that $K$ has \emph{IVP} (the \emph{Intermediate Value Property}) if for all
$P(Y)\in K\{Y\}$ and $f< g$ in $K$ with $P(f) <0 <  P(g)$ there is a~$\phi\in K$ such that $f < \phi < g$ and $P(\phi)=0$. Restricting this to $P$ of order~$\le r$, where~$r\in \N$, gives the notion of $r$-IVP. Thus
$K$ having $0$-IVP is equivalent to~$K$ being real closed as an
ordered field. In particular, if $K$ has $0$-IVP, then the
$H$-asymptotic couple~$(\Gamma, \psi)$ of $K$ is divisible. 
From \cite[Section~2.4]{ADH} recall our convention that $K^>=\{a\in K:a>0\}$, and similarly with $<$ replacing $>$.

\begin{lemma}\label{ivp1} Suppose $\Gamma\ne \{0\}$ and $K$ has 
$1$-$\operatorname{IVP}$. Then
$\der K=K$, $(K^{>})^\dagger=(K^{<})^\dagger$ is a convex subgroup of $K$, $\Psi:=\big\{\psi(\gamma):\gamma\in\Gamma^{\ne}\big\}$ has no largest element, and $\Psi$ is convex in $\Gamma$.
\end{lemma}
\begin{proof} We have $y'=0$ for $y=0$, and $y'$ takes arbitrarily large positive values in~$K$ as $y$ ranges over
$K^{>\mathcal{O}}=\{a\in K:a>\mathcal O\}$, since by \cite[Lemma~9.2.6]{ADH} the set
$(\Gamma^{<})'$ is coinitial in~$\Gamma$. Hence
$y'$ takes all positive values on $K^{>}$, and therefore also all negative values on~$K^{<}$. Thus $\der K=K$. Next, let $a,b\in K^{>}$, and suppose
$s\in K$ lies strictly between
$a^\dagger$ and $b^\dagger$. Then $s=y^\dagger$ for some $y\in K^{>}$
strictly between $a$ and $b$; this follows by noting that for
$y=a$ and $y=b$ the signs of
$sy-y'$ are opposite.

Let $\beta\in\Psi$ 
and take $a\in K$ with $v(a')=\beta$. Then $a\succ 1$, since
$a\preceq 1$ would give $v(a')>\Psi$. Hence for $\alpha=va<0$
we have $\alpha+\alpha^\dagger=\beta$, so $\alpha^\dagger>\beta$.
Thus $\Psi$ has no largest element. 
Therefore the set $\Psi$ is convex in $\Gamma$.  \end{proof}  

\noindent
Thus the ordered differential field $\T_{\log}$ of   logarithmic transseries,
\cite[Appendix~A]{ADH}, does not have $1$-IVP, although it is a newtonian 
$\upo$-free $H$-field.

Does IVP imply that $K$ is an $H$-field? No: take an $\aleph_0$-saturated elementary
extension of $\T$ and let $\Delta$ be as in \cite[Example~10.1.7]{ADH}. Then
the $\Delta$-coarsening of $K$ is a pre-$H$-field with IVP and nontrivial value group, and has a gap, but it is not an $H$-field.
On the other hand:

\begin{lemma}\label{ivp2} Suppose $K$ has $1$-$\operatorname{IVP}$ and has no gap. Then $K$ is an $H$-field.
\end{lemma}
\begin{proof} In \cite[Section~11.8]{ADH} we defined
$$\text{I}(K):= \{y\in K:\ \text{$y\preceq f'$ for some $f\in \mathcal{O}$}\}.$$ Since $K$ has no gap, we have
$$\der\smallo\ \subseteq\ \text{I}(K)\ =\ \{y\in K:\ \text{$y\preceq f'$ for some $f\in \smallo$}\}.$$
Also $\Gamma\ne \{0\}$, and so $(\Gamma,\psi)$ has asymptotic integration by Lemma~\ref{ivp1}. We show that $K$ is an $H$-field
by proving $\text{I}(K) =\der\smallo$, so let $g\in \text{I}(K)$, $g<0$. Since $(\Gamma^{>})'$ has no least element we can take
positive $f\in \smallo$ such that $f'\succ g$. Since
$f'<0$, this gives $f' < g$. Since $(\Gamma^{>})'$ is cofinal in 
$\Gamma$ we can also take positive $h\in \smallo$ such that $h'\prec g$, which in view of $h'<0$ gives $g < h'$.
Thus $f' < g < h'$, and so $1$-IVP yields $a\in \smallo$ with $g=a'$. 
\end{proof}

\noindent
We refer to Sections~11.6 and~14.2 of \cite{ADH} for the definitions of {\it $\upl$-freeness}\/ and {\it $r$-newtonianity}\/ ($r\in\N$). From the introduction we recall that
$\omega(z):=-2z'-z^2$.

\begin{lemma}\label{ivp3} Suppose $K$ is an $H$-field, $\Gamma\ne \{0\}$, and $K$ has 
$1$-$\operatorname{IVP}$. Then $K$ is $\upl$-free and $1$-newtonian, and the subset $\omega(K)$ of $K$ is downward closed.  
\end{lemma}
\begin{proof}  First we note that $K$ has (asymptotic) integration, by Lemma~\ref{ivp1}. Assume towards a contradiction that $K$ is
not $\upl$-free. We can arrange that~$K$ has small derivation, and thus $K$ has an element $x\succ 1$ with $x'=1$, and so $x> C$.
This leads to a pc-sequence $(\upl_\rho)$ and an element
$s\in K$ such that $\upl_{\rho} \leadsto -s$ with $\upl_{\rho}\sim x^{-1}$ for all $\rho$. Hence  $s\sim -x^{-1}$, and
$s$ creates a gap
over $K$ by \cite[Lemma~11.5.14]{ADH}. Now note that for
$P:= Y'+sY$ we have $P(0)=0$ and $P(x^2)=2x+sx^2\sim x$, so
by $1$-IVP we have $P(y)=1$ for some $y\in K$, contradicting
\cite[Lemma~11.5.12]{ADH}. 

Let $P\in K\{Y\}$ of order at most~$1$ have Newton degree~$1$; we have to show that $P$ has a zero in $\mathcal{O}$. We know that $K$ is $\upl$-free, so by \cite[Proposition~13.3.6]{ADH} we
can pass to an elementary extension, 
compositionally conjugate, and divide by an element of $K^\times$ to arrange that $K$ has small derivation and
$P=D+R$ where $D=cY+d$ or $D=cY'$ with $c,d\in C$, $c\ne 0$, and where $R\prec^{\flat} 1$. Then $R(a)\prec^{\flat} 1$ for all $a\in \mathcal{O}$. If $D= cY+d$, then we can take $a,b\in C$ with $D(a)<0$ and $D(b)>0$, which in view of $R(a)\prec D(a)$ and
$R(b) \prec D(b)$ gives $P(a)<0$ and $P(b)>0$, and so
$P$ has a zero strictly between $a$ and $b$, and thus a zero in $\mathcal{O}$. Next, suppose $D=cY'$. Then we take $t\in \smallo^{\ne}$ with $v(t^\dagger)=v(t)$, that is, $t'\asymp t^2$, so
$$P(t)\ =\ ct'+R(t), \quad P(-t)\ =\ -ct'+ R(-t), \qquad 
R(t),\ R(-t)\ \prec\ t'.$$
Hence $P(t)$ and $P(-t)$ have opposite signs, so
$P$ has a zero strictly between $t$ and~$-t$, and thus $P$
has a zero in $\mathcal{O}$.

From $\omega(z)=-z^2-2z'$ we see that $\omega(z) \to -\infty$ as
$z\to +\infty$ and as $z\to -\infty$ in~$K$, so $\omega(K)$ is downward closed by $1$-IVP.  
\end{proof} 

\noindent
For results involving $2$-IVP we need a minor variant of \cite[Lemma~11.8.31]{ADH}. Here
$\Upg(K) = \{ a^\dagger: a\in K\setminus\mathcal O\}$ as in~\cite[Section~11.8]{ADH},
and the superscripts~$\uparrow$,~$\downarrow$ indicate upward, respectively downward, closure, as in~\cite[Section~2.1]{ADH}.

\begin{lemma}\label{ivp4-} Let $K$ be an $H$-field with asymptotic integration. Then 
$$K^{>}\ =\ \operatorname{I}(K)^{>}\cup\Upg(K)^{\uparrow}, \qquad \sigma\big(K^{>}\setminus \Upg(K)^{\uparrow}\big)\ \subseteq\ \omega(K)^{\downarrow}.$$ 
\end{lemma}
\begin{proof} If $a\in K$, $a > \text{I}(K)$, then $a\ge b^\dagger$ for some $b\in K^{\succ 1}$, and thus $a\in \Upg(K)^{\uparrow}$. 
The inclusion involving $\sigma$ now follows as in the proof of 
\cite[Lemma 11.8.13]{ADH}. 
\end{proof} 

\noindent
The concept of \emph{$\upo$-freeness}\/ is introduced in \cite[Section~11.7]{ADH}.

\begin{lemma}\label{ivp4} Suppose $K$ is an $H$-field, $\Gamma\ne \{0\}$, and $K$ has
$2$-$\operatorname{IVP}$. Then $K$ is $2$-newtonian,
the operator $\der^2-a$ splits over $K[\imag]$ for all $a\in K$,
and $K$ is $\upo$-free. 
\end{lemma}
\begin{proof} Let $P\in K\{Y\}$ of order at most $2$ have Newton degree $1$; we have to show that $P$ has a zero in $\mathcal{O}$. Lemma~\ref{ivp3} tells us that $K$ is $\upl$-free, and 
in view of \cite[Corollary~13.3.7]{ADH} and $2$-IVP this allows us to
repeat the argument in the proof of that lemma for differential polynomials of order at most $1$ so that it applies to our~$P$ of order at most $2$. Thus $K$ is $2$-newtonian.

By \cite[Section~5.2]{ADH} it remains to show that
 $K=\omega(K)\cup \sigma(K^\times)$. In view of Lemma~\ref{ivp1} we can arrange by compositional conjugation
that $a^\dagger=-1$ for some~${a\prec 1}$ in $K^{>}$. Below we fix such $a$. Let $f\in K$; our job is to show that $f\in \omega(K)\cup \sigma(K^\times)$.
Since $\omega(0)=0$, we do have $f\in \omega(K)$ if $f\le 0$,
by Lemma~\ref{ivp3}. So assume $f> 0$;
we show that then $f\in \sigma(K^{>})$. Now  for $y\in K^{>}$, $f=\sigma(y)$ is equivalent (by multiplying with $y^2$) to $P(y)=0$, where $$P(Y)\ :=\ 2YY'' -3(Y')^2+Y^4-fY^2\in K\{Y\}.$$
See also \cite[Section~13.7]{ADH}. We have $P(0)=0$ and $P(y) \to +\infty$ as $y\to +\infty$ (because of the term $y^4$). In view of $2$-IVP it will suffice to show that
for some $y>0$ in $K$ we have $P(y)<0$. Now with
$y\in K^{>}$ and $z:=-y^\dagger$ we have
$$P(y)=y^2\big(\sigma(y)-f\big)=y^2\big(\omega(z)+y^2-f\big),$$
hence 
$$P(a)=a^2\big(\omega(1)+a^2-f\big)=a^2(-1+a^2-f)<0.$$ 
As to $\upo$-freeness, this now follows from Lemma~\ref{ivp4-} and \cite[Corollary~11.8.30]{ADH}.
\end{proof} 

\noindent
It follows that Liouville closed $H$-fields having $2$-$\operatorname{IVP}$ are Schwarz closed
as defined in~\cite[Sec\-tion~11.8]{ADH}. (There exist  $H$-fields with a non-trivial derivation that have~$\operatorname{IVP}$ but are not Liouville closed; see \cite[Section~14]{AvdD}.) 

\begin{cor}\label{LIVP} Suppose $K$ is an $H$-field, $\Gamma\ne \{0\}$, and $K$ has 
$\operatorname{IVP}$. Then $K$ is $\upo$-free and newtonian. 
\end{cor} 
\begin{proof} Showing that every $P\in K\{Y\}$ of Newton degree $1$ has a zero in $\mathcal{O}$ is done just as in the proof of Lemma~\ref{ivp3}.  
\end{proof}

\begin{cor} Let $K$ be a Liouville closed $H$-field. Then 
$$ \text{$K$   has $\operatorname{IVP}$}\ \Longleftrightarrow\ \text{$K$  is $\upo$-free and newtonian.}$$
\end{cor}
\begin{proof} The forward direction is part of Corollary~\ref{LIVP}. For the backward direction we appeal to the main results
from the book~\cite{JvdH}  to the effect that $\T_{\text{g}}$, the ordered differential field of grid-based transseries (cf.~\cite[Appendix~A]{ADH}), is a
newtonian Liouville closed $H$-field with small derivation, and has IVP. In particular, it is a model of the theory $T^{\text{nl}}_{\text{small}}$, which we mentioned in the introduction.  This theory is complete by~\cite[Corollary~16.6.3]{ADH}, 
so every model of it has IVP. If $K$ is $\upo$-free and newtonian but its derivation is not small, then it nevertheless  has
IVP: some compositional conjugate $K^{\phi}$ with $\phi \in K^{>}$ has small derivation and is Liouville closed, $\upo$-free and newtonian.
\end{proof}

\section{Preliminaries on Hardy Fields}\label{bh}

\noindent
We begin with some results from Boshernitzan~\cite{B} on ordered fields of germs of continuous functions. Next we prove some easy facts about extending ordered fields inside an ambient partially
ordered ring, as needed later.

\subsection*{Germs of continuous functions} As in \cite[Section 9.1]{ADH} we let $\mathcal{G}$ be the ring
of germs at $+\infty$ of real-valued functions whose domain is
a subset of $\R$
containing an interval $(a, +\infty)$, $a\in \R$; the domain may vary 
and the ring operations are defined as usual. 
If $g\in \mathcal{G}$ is the germ of a real-valued function 
on a subset of $\R$ containing an interval $(a, +\infty)$, $a\in \R$, then we simplify notation
by letting $g$ also denote this function if the resulting ambiguity is harmless.
With this convention, given a property~$P$ of real numbers
and $g\in \mathcal{G}$ we say that {\em $P\big(g(t)\big)$ holds eventually\/} if~$P\big(g(t)\big)$ holds for all sufficiently large real $t$. 
We identify each real number $r$ with the germ at~$+\infty$ of the function 
$\R\to \R$ that takes the constant value $r$. This  makes the field $\R$ into a subring of $\mathcal{G}$. We call a germ $g\in \mathcal{G}$ {\em continuous} if it is the germ of
a continuous function $(a,+\infty)\to \R$ for some $a\in \R$, and we let $\mathcal{C}\supseteq \R$ be the subring of $\mathcal{G}$ consisting of the continuous germs $g\in \mathcal{C}$. We let $x$ denote the germ at $+\infty$ of the identity function on $\R$.

\subsection*{Asymptotic relations on $\C$.} Note that the multiplicative group $\C^\times$ of $\C$ consists of the $f\in \C$ such that $f(t)\ne 0$, eventually. Thus for $f\in \C^\times$, either $f(t)>0$, eventually, or $f(t)<0$, eventually. Although $\C$ is not a valued field, it will be convenient to equip $\C$ with the asymptotic relations $\preceq$,~$\prec$,~$\sim$ (which are defined on any valued field) as follows: for $f,g\in \C$,
\begin{align*} f\preceq g\ &:\Longleftrightarrow\ \text{there exists $c\in \R^{>}$ such that eventually}\  |f(t)|\le c|g(t)|,\\
f\prec g\ &:\Longleftrightarrow\ \text{$g\in \C^\times$ and $\lim_{t\to \infty} f(t)/g(t)=0$},\\
f\sim g\ &:\Longleftrightarrow\ \text{$g\in \C^\times$ and
$\lim_{t\to \infty} f(t)/g(t)=1$}\\ 
&\ \Longleftrightarrow\ f-g\prec g.
\end{align*}
Thus $\preceq$ is a transitive and reflexive binary relation on
$\C$, and $\sim$ is an equivalence relation on $\C^\times$.
Moreover, for $f,g,h\in \C$ we have
$$f\prec g\ \Rightarrow\ f\preceq g, \qquad f\preceq g \prec h\ \Rightarrow\ f\prec h, \qquad f\prec g \preceq h\ \Rightarrow\ f\prec h.$$ 
Note that $\prec$ is a transitive binary relation
on $\C$.  
For $f,g\in \C$ we also set
$$f\asymp g:\ \Leftrightarrow\ f\preceq g \text{ and }g\preceq f,\qquad f\succeq g:\ \Leftrightarrow\ g\preceq f,\qquad f\succ g:\ \Leftrightarrow\ g\prec f,
$$
so $\asymp$ is an equivalence relation on $\C$.

\subsection*{Subfields of $\mathcal{C}$}
Let $K$ be a subfield of $\mathcal{C}$, that is, a subring
of $\mathcal{C}$ that happens to be a field. (In the introduction we called such $K$ a {\em Hausdorff field}\/.) Then $K$ itself has the subfield $K\cap \R$. Every nonzero $f\in K$ has a multiplicative inverse in $K$, so eventually $f(t)\neq 0$, hence
either eventually $f(t)<0$ or eventually $f(t)>0$ (by eventual continuity of $f$).
We make~$K$ an ordered field by declaring 
$$f>0\ :\Longleftrightarrow\ \text{$f(t) > 0$, eventually.}$$
We now have \cite[Propositions 3.4 and 3.6]{B}: 

\begin{lemma}\label{b1} Let $K^{\operatorname{rc}}$ consist of the
$y\in \mathcal{C}$ with $P(y)=0$ for some $P(Y)\in K[Y]^{\ne}$. Then $K^{\operatorname{rc}}$ is the
unique real closed subfield of $\mathcal{C}$ that
extends $K$ and is algebraic over~$K$. In particular, $K^{\operatorname{rc}}$ is a real closure of the ordered field $K$.  
\end{lemma} 

\noindent
In \cite{B} this lemma assumes $K\supseteq \R$, but this is not really needed in the proof. 
The ordered field $K$ has a convex subring 
$$ \mathcal{O}\ =\ 
\big\{f\in K :\ \text{$|f|\le n$ for some $n$}\big\},$$
which is a valuation ring of $K$, and we
consider $K$ accordingly as a valued ordered field. Restricting
$\preceq$,~$\prec$,~$\sim$ from the previous subsection to $K$ gives exactly the asymptotic relations $\preceq$,~$\prec$,~$\sim$ on $K$
that it comes equipped with as a valued field.

\subsection*{Composition and compositional inversion.} Let $g\in \C$ be eventually strictly increasing with $\lim_{t\to +\infty} g(t)=+\infty$. Then its compositional inverse $g^{\inv}\in \C$ is given by $g^{\inv}\big(g(t)\big)=t$, eventually, and the composition operation
$$f\mapsto f\circ g\ :\ \C \to \C, \qquad (f\circ g)(t)\ :=\  f\big(g(t)\big)\ \text{ eventually},$$
is an automorphism of the ring $\C$
that is the identity on the subring $\R$, with inverse~$f\mapsto f\circ g^{\inv}$.  In particular, $g\circ g^{\inv}=g^{\inv}\circ g=x$, and $f\mapsto f\circ g$ maps each
subfield $K$ of $\C$ isomorphically (as an ordered field)
onto the subfield $K\circ g$ of $\C$. Note that if the subfield $K$ of $\C$ contains $x$, then $K\circ g$ contains $g$.

\subsection*{Extending ordered fields inside an ambient partially ordered ring}
Let $R$ be a commutative ring with $1\ne 0$, equipped with a translation-invariant
partial ordering $\le$ such that $r^2\ge 0$ for all $r\in R$, and
$rs\ge 0$
for all $r,s\in R$ with $r,s\ge 0$. It follows that for $a,b,r\in R$ we have: if $a \le b$ and $r\ge 0$,
then $ar\le br$; if $a$ is a unit and $a>0$, then $a^{-1}=a\cdot (a^{-1})^2>0$; if $a,b$ are units, and $0 < a \le b$, then $0 < b^{-1}\le a^{-1}$. 
Relevant cases: 
$R=\mathcal{G}$ and $R=\mathcal{C}$, with partial ordering given by $$f\le g\ :\ \Longleftrightarrow\ f(t) \le g(t),  \text{ eventually}.$$ 
Call a subset $K$ of $R$ {\em totally ordered\/} if the partial ordering of $R$ induces a total ordering on $K$.
An {\em ordered subfield of $R$} is a subfield $K$
of $R$ that is totally ordered as a subset of $R$; note that then $K$ equipped with the induced partial ordering
is indeed an ordered field, in the usual sense of that term.
(Thus any subfield of $\mathcal{C}$ with the above partial ordering is an ordered subfield of $\mathcal{C}$.) 

We identify $\Z$ with its image in $R$ via the unique ring embedding $\Z \to R$, and this makes $\Z$ with its usual ordering into an ordered subring of $R$. 

\begin{lemma}\label{pr0} Assume $D$ is a totally ordered
subring of $R$ and every nonzero element of $D$ is a unit of $R$.
Then $D$ generates an ordered subfield $\Frac{D}$ of $R$.
\end{lemma}
\begin{proof} It is clear that $D$ generates a subfield $\Frac{D}$ of $R$. For $a\in D$, $a>0$, we have $a^{-1}>0$. It follows that $\Frac{D}$ is totally ordered.
\end{proof}

\noindent
Thus if every $n\ge 1$ is a unit of $R$, then we may identify $\Q$ with its image in $R$ via the unique
ring embedding $\Q\to R$, making $\Q$ into an ordered subfield of $R$.

\begin{lemma}\label{pr2} Suppose $K$ is an ordered subfield of $R$,
all $g\in R$ with $g>K$ are units of $R$, and $K<f\in R$. Then we have an ordered subfield $K(f)$ of $R$. 
\end{lemma}
\begin{proof} 
For $P(Y)\in K[Y]\setminus K$ of degree $d\ge 1$
with leading coefficient $a>0$ we have $P(f)=af^d(1+\epsilon)$
with $-1/n < \epsilon < 1/n$ for all $n\ge 1$, in particular,
$P(f)>K$ is a unit of $R$. It remains to appeal to Lemma~\ref{pr0}.
\end{proof}

\begin{lemma}\label{pr1} Assume $K$ is a real closed ordered subfield of $R$.
Let $A$ be a nonempty downward closed subset of $K$ such that
$A$ has no largest element and $B:=K\setminus A$ is nonempty and
has no least element. Let $f\in R$ be such that $A<f<B$. Then the subring $K[f]$ has the following properties: \begin{enumerate}
\item[$\rm{(i)}$] $K[f]$ is a domain;
\item[$\rm{(ii)}$] $K[f]$ is totally ordered;
\item[$\rm{(iii)}$] $K$ is cofinal in $K[f]$;
\item[$\rm{(iv)}$] for all $g\in K[f]\setminus K$ and $a\in K$, if $a<g$, then $a < b<g$ for some $b\in K$, and if $g<a$, then 
$g<b<a$ for some $b\in K$. 
\end{enumerate}
\end{lemma} 
\begin{proof} Let
$P\in K[Y]\setminus K$; to obtain (i) and (ii) it suffices to show that then
$P(f) < 0$ or $P(f)>0$.  We have $$P(Y)\ =\ cQ(Y)(Y-a_1)\cdots(Y-a_n)$$ where
$c\in K^{\ne}$,  $Q(Y)$  is  a product of monic quadratic irreducibles in $K[Y]$, and  $a_1,\dots, a_n\in K$. This gives $\delta\in K^{>}$ such that $Q(r)\ge \delta$ for all $r\in R$.  Assume~$c>0$. (The case $c<0$ is handled similarly.) We can arrange that $m\le n$ is such that
$a_i\in A$ for $1\le i \le m$ and $a_j\in B$ for $m < j\le n$.
Take $\epsilon>0$ in $K$ such that
$a_i + \epsilon \le f$ for $1\le i\le m$ and 
$f\le a_j-\epsilon$ for $m < j\le n$.
Then $$P(f)\ =\ cQ(f)(f-a_1)\cdots (f-a_m)(f-a_{m+1})\cdots(f-a_n),$$
and $(f-a_1)\cdots(f-a_m) \ge \epsilon^m$. If
$n-m$ is even, then $(f-a_{m+1})\cdots(f-a_n)\ge \epsilon^{n-m}$, so $P(f)\ge a\delta\epsilon^n >0$. If $n-m$ is odd, then
$(f-a_{m+1})\cdots(f-a_n)\le -\epsilon^{n-m}$,
so $P(f) \le -a\delta\epsilon^n < 0$. These estimates also yield (iii) and (iv). 
\end{proof}

\begin{lemma}\label{pr3} With $K$,~$A$,~$f$ as in Lemma~\ref{pr1}, suppose all $g\in R$ with $g\ge 1$ are units of $R$. Then we have an ordered subfield $K(f)$ of $R$ such that {\rm(iii)} and {\rm(iv)} of Lemma~\ref{pr1} go through for $K(f)$ in place of $K[f]$.
\end{lemma} 
\begin{proof} 
Note that if $g\in R$ and $g\ge \delta\in K^{>}$, then 
$g\delta^{-1}\ge 1$, so $g$ is a unit of $R$ and
$0<g^{-1}\le \delta^{-1}$. For   
$Q\in K[Y]^{\ne}$ with $Q(f)>0$ we can take $\delta\in K^{>}$ such that $Q(f)\ge \delta$, and thus $Q(f)$ is a unit
of $R$ and $0 <Q(f)^{-1}\le \delta^{-1}$. Thus we have an
ordered subfield
$K(f)$ of $R$ by Lemma~\ref{pr0}, and 
the rest now follows easily.  
\end{proof}

\subsection*{Adjoining pseudolimits and increasing the value group} Let $K$ be a real closed subfield of $\mathcal{C}$, and  view $K$ as an ordered valued field as
before. Let $(a_{\rho})$ be a strictly increasing divergent pc-sequence in $K$.
Set 
$$A\ :=\ \{a\in K:\ \text{$a< a_{\rho}$ for some $\rho$}\}, \qquad
B\ :=\ \{b\in K:\ \text{$b>a_{\rho}$ for all $\rho$}\},$$ so
$A$ is nonempty and downward closed without a largest element.
Moreover, $B=K\setminus A$ is nonempty and has no least element, since a
least element of $B$ would be a limit and thus a pseudolimit of
$(a_{\rho})$. Let $f\in \mathcal{C}$ satisfy $A<f<B$.
Then we have an ordered subfield $K(f)$ of $\mathcal{C}$, and:

\begin{lemma} \label{ps1} $K(f)$ is an immediate valued field extension of $K$ with $a_{\rho} \leadsto f$.
\end{lemma}
\begin{proof} We can assume that $v(a_{\tau}-a_{\sigma})>v(a_{\sigma}-a_{\rho})$ for all indices $\tau>\sigma>\rho$. Set $d_{\rho}:= a_{s(\rho)}-a_{\rho}$ ($s(\rho):=$~successor of $\rho$).
Then $a_{\rho}+2d_{\rho}\in B$ for all indices~$\rho$; see the discussion
preceding~\cite[Lem\-ma~2.4.2]{ADH}. It then follows from that lemma that $a_{\rho} \leadsto f$. Now $(a_{\rho})$ is a divergent pc-sequence in the henselian valued field $K$, so it is of transcendental type over $K$, and thus $K(f)$ is an immediate
extension of~$K$.  
\end{proof}

\begin{lemma}\label{ps2} Suppose $K$ is a subfield of 
$\mathcal{C}$ with divisible
value group $\Gamma=v(K^\times)$. Let $P$ be a nonempty upward closed subset of $\Gamma$, and let $f\in \mathcal{C}$ be such that 
$a < f$ for all $a\in K^{>}$ with 
$va\in P$, and $f<b$ for all $b\in K^{>}$ with $vb < P$.
Then $f$ generates a subfield $K(f)$ of $\mathcal{C}$, with 
$P >vf > Q$, $Q:= \Gamma\setminus P$.  
\end{lemma}
\begin{proof} For any positive $a\in K^{\text{rc}}$ there is
$b\in K^{>}$ with $a\asymp b$ and $a< b$, and also an element
$b\in K^{>}$ with $a\asymp b$ and $a>b$. Thus we can replace $K$ by $K^{\text{rc}}$ and arrange in this way that $K$ is real closed.
 Set 
$$A\ :=\ \{a\in K:\ a\le 0 \text{ or }va\in P\}, \qquad B:= K\setminus A.$$
Then we are in the situation of Lemma~\ref{pr1} for $R=\mathcal{C}$,   so by that lemma and Lemma~\ref{pr3} we have an ordered subfield
$K(f)$ of $\mathcal{C}$. Clearly then $P> vf >Q$.     
\end{proof}

\subsection*{Notational conventions on functions and germs} 
Let $r$ range over $\N\cup\{\infty\}$, and let $U$ be a nonempty open subset of $\R$. 
Then $\C^r(U)$ denotes the $\R$-algebra of $r$-times
continuously differentiable functions $U\to \R$, with the usual pointwise defined algebra operations. (We use ``$\C$'' instead of ``$C$'' since $C$ will often denote the constant field of a
differential field.) For $r=0$ this is the $\R$-algebra~$\C(U)$ of continuous real-valued functions on $U$, so 
$$\C(U)\ =\ \C^0(U)\ \supseteq\ \C^1(U)\ \supseteq\ \C^2(U)\ \supseteq\ 
\cdots\ \supseteq\ \C^{\infty}(U).$$ For $r\ge 1$ we have the derivation $f\mapsto f'\colon \C^r(U)\to \C^{r-1}(U)$ (with $\infty-1:=\infty$). This makes $\C^{\infty}(U)$ a differential ring, with
its subalgebra
$\C^{\omega}(U)$ 
of real-analytic functions $U\to \R$ as a differential subring. The algebra operations
on the algebras below are also defined pointwise.   

Let $a$ range over $\R$. 
Then $\Car$ denotes the $\R$-algebra of functions $[a,+\infty) \to \R$ that extend to a function in $\C^r(U)$ for some open $U\supseteq [a,+\infty)$. Thus $\Caz$ is the $\R$-algebra of real-valued continuous functions on $[a,+\infty)$, and
$$\Caz\ \supseteq\ \Cao\ \supseteq\ \Cat\ \supseteq\ \cdots\ \supseteq \Cainf.$$  
We also have the subalgebra $\Caom$ of $\Cainf$, consisting of the functions
$[a,+\infty)\to \R$ that extend to a real-analytic function
$U \to \R$ for some open $U\supseteq [a,+\infty)$. For~$r\ge 1$ we have the derivation $f\mapsto f'\colon \C^r_a\to \C^{r-1}_a$. This makes $\C^{\infty}_a$ a differential ring with~$\C^{\omega}_a$  as a differential subring.

For each of the algebras $A$ above we also consider its complexification $A[\imag]$ which consists by definition of the
$\Co$-valued functions $f=g+h\imag$ with $g,h\in A$, so
$g=\operatorname{Re} f$ and $h= \operatorname{Im} f$ for such $f$. We consider $A[\imag]$ as a $\Co$-algebra with respect to the natural pointwise defined algebra operations. We identify each complex number with
the corresponding constant function to make $\Co$ a subfield of
$A[\imag]$ and $\R$ a subfield of~$A$. (This justifies the notation $A[\imag]$.) For $r\ge 1$ we extend $g\mapsto g'\colon \C^r_a\to \C^{r-1}_a$ to the derivation 
$$g+h\imag\mapsto g'+h'\imag\ :\  \C^r_a[\imag] \to 
\C^{r-1}_a[\imag] \qquad (g,h\in \C^r_a[\imag]),$$
which for $r=\infty$ makes $\C^{\infty}_a$ a differential subring of $\C^{\infty}_a[\imag]$.  
We also use the map
$$f \mapsto f^\dagger:=f'/f\ \colon \ \C^1_a[\imag]^\times=\big(\C^1_a[\imag]\big)\!^\times \to \C^{0}_a[\imag],$$
with 
$$(fg)^\dagger=f^\dagger+g^\dagger\qquad\text{ for $f,g\in \C^1_a[\imag]^\times$,}$$
in particular the fact that $f\in \C^1_a[\imag]^\times$ 
and $f^\dagger\in \C^0_a[\imag]$ are related by
$$ f(t)\ =\ f(a)\exp\left[\int_a^t f^\dagger(s)\,ds\right] \qquad (t\ge a).$$ 
Let $\Gr$ be the partially ordered subring of 
$\C$ consisting of the germs at $+\infty$ of the functions in $\bigcup_a \Car$; thus $\Gz=\C$ consists of the germs at $+\infty$ of the continuous real valued functions on intervals $[a,+\infty)$, $a\in \R$. Note that $\Gr$ with its partial ordering satisfies the conditions on $R$ from the previous subsection. Also, every~$g\ge 1$ in~$\Gr$ is a unit of $\Gr$, so
Lemmas~\ref{pr2} and \ref{pr3} apply to ordered subfields of $\Gr$. We have 
$$\Gz\ \supseteq\ \Go\ \supseteq\ \Gt\ \supseteq\ \cdots\ \supseteq\ \Ginf,$$
and we set $\Gi:= \bigcap_r \Gr$. Thus $\Gi$ is naturally a
differential ring with $\R$ as its ring of constants. Note that
$\Gi$ has~$\Ginf$ as a differential subring. The differential ring $\Ginf$ has in turn
the differential subring~$\Gom$, whose elements are the
germs at~$+\infty$ of the functions in $\bigcup_a \Caom$.

\subsection*{Second-order differential equations} Let $f\in \C^0_a$, that is, $f\colon [a,\infty) \to \R$ is continuous. We consider the differential equation
$$ Y'' + fY\ =\ 0. $$
The solutions~${y\in \C^2_a}$ form
an $\R$-linear subspace $\Sol(f)$ of $\C^2_a$. The
solutions~${y\in \C^2_a[\imag]}$ are the $y_1+y_2\imag$ with
$y_1, y_2\in \Sol(f)$ and form
a $\Co$-linear subspace $\Sol_{\Co}(f)$ of~$\C^2_a[\imag]$. For any complex numbers
$c$,~$d$ there is a unique solution $y\in \C^2_a[\imag]$ with~${y(a)=c}$ and $y'(a)=d$, and the map that assigns to $(c,d)\in \Co^2$
this unique solution is an isomorphism $\Co^2\to \Sol_{\Co}(f)$ of 
$\Co$-linear spaces; it restricts to an $\R$-linear bijection~$\R^2\to \Sol(f)$. Induction on $r\in \N$ shows: 
$f\in \C^r_a\Rightarrow \Sol(f)\subseteq \C^{r+2}_a$. Thus $f\in \C^\infty_a\Rightarrow \Sol(f)\subseteq \C^\infty_a$. It is also well-known
that $f\in \C^\omega_a\Rightarrow \Sol(f)\subseteq \C^\omega_a$. 
From~\cite[Chapter~2, Lemma~1]{Bellman} we recall:

\begin{lemma}[Gronwall's Lemma] 
Let the constant $C\in\R^{\ge}$ and the functions $v,y\in\mathcal C^0_a$ be such that $v(t),y(t)\geq 0$ for all $t\geq a$ and
$$y(t)\ \leq\ C+ \int_a^t v(s)y(s)\,ds\quad\text{for all $t\geq a$.}$$
Then
$$y(t)\ \leq\ C\exp\left[ \int_a^t v(s)\,ds\right]\quad\text{for all $t\geq a$.}$$
\end{lemma}

\noindent
{\it In the rest of this subsection we assume that
$a\geq 1$ and that $c\in\R^>$  is such that $|f(t)| \leq c/t^2$ for all $t\geq a$.}
Under this hypothesis, the lemma above yields the following bound on the growth of the solutions $y\in\Sol(f)$; the proof we give is similar to that of~\cite[Chap\-ter~6, Theorem~5]{Bellman}.

\begin{prop}\label{prop:bound} Let $y\in\Sol(f)$. Then there is $C\in\R^\geq$ such that $|y(t)| \leq Ct^{c+1}$ 
and $|y'(t)|\leq Ct^c$
for all $t\geq a$.
\end{prop}
\begin{proof}
Let $t$ range over $[a,+\infty)$.
Integrating $y''=-fy$ twice between $a$ and $t$, we obtain constants $c_1$, $c_2$ such that
for all $t$,
$$y(t)\ =\ c_1+c_2t -\int_a^t\int_a^{t_1} f(t_2)y(t_2)\,dt_2\,dt_1\ =\ c_1+c_2t-\int_a^t (t-s)f(s)y(s)\, ds$$
and hence, with $C:=|c_1|+|c_2|$,
$$|y(t)|\ \leq\ Ct + t\int_a^t |f(s)|\cdot|y(s)|\, ds,$$
so
$$\frac{|y(t)|}{t}\ \leq\ C + \int_a^t s|f(s)|\cdot\frac{|y(s)|}{s}\, ds.$$
Hence by the lemma above, 
$$\frac{|y(t)|}{t}\ \leq\ C\exp\left[\int_a^t s|f(s)|\,ds\right]\ \leq\
C\exp\left[\int_1^t c/s\,ds\right]\ =\ Ct^c$$
and thus
$|y(t)|\leq Ct^{c+1}$. Now  
$$y'(t)=c_2-\int_a^t f(s)y(s)\,ds$$
and thus
$$|y'(t)|\ \leq\ |c_2|+\int_a^t |f(s)y(s)|\,ds\ \leq\ C+Cc\int_1^t s^{c-1}\,ds\ =\
C+Cc\left[\frac{t^c}{c}-\frac{1}{c}\right]\ =\ Ct^c.
$$
\end{proof} 

\noindent
Let $y_1, y_2\in \Sol(f)$ be $\R$-linearly independent. The Wronskian $w:= y_1y_2'-y_1'y_2$ satisfies
$w'=0$ (Abel's identity), so $w\in \R^\times$. It follows that $y_1$ and $y_2$ cannot be simultaneously very small:

\begin{lemma}\label{lem:bound} There is a positive constant $d$ such that
$$\max\big(|y_1(t)|, |y_2(t)|\big)\ \ge\ dt^{-c} \quad \text{ for all $t\ge a$.}$$
\end{lemma}
\begin{proof} 
Proposition~\ref{prop:bound} yields $C\in \R^>$ such that  $|y_i'(t)|\le Ct^c$ for
$i=1,2$ and all~$t\ge a$. Hence $|w|\le 2\max\big(|y_1(t)|, |y_2(t)|\big)Ct^c$ for $t\ge a$, so
\equationqed{ \max\big(|y_1(t)|, |y_2(t)|\big)\ \ge\ \frac{|w|}{2C}t^{-c}\qquad(t\ge a).}
\end{proof}

\begin{cor}\label{cor:bound}
Set $y:=y_1+y_2\imag$ and $z:=y^\dagger$. Then for some $D\in\R^>$,  
$$|z(t)|\ \leq\ Dt^{2c}\quad\text{ for all $t\geq a$.}$$
\end{cor}
\begin{proof}
Take $C$ as in the proof of Lemma~\ref{lem:bound}, and $d$ as in that lemma. Then
$$|z(t)|\ =\ \frac{|y_1'(t)+y_2'(t)\imag|}{|y_1(t)+y_2(t)\imag|}\ \leq\  \frac{|y_1'(t)|+|y_2'(t)|}{\max\big(|y_1(t)|, |y_2(t)|\big)}\ \leq\ \left(\frac{2C}{d}\right)t^{2c}$$
for $t\geq a$.
\end{proof}

\subsection*{Changing variables} Let now $K$ be a differential field, $f\in K$, and consider the
differential polynomial $P(Y):= 4Y''+fY$. (The factor $4$ is to simplify certain expressions, in conformity with \cite[Section~9.2]{ADH}.) Which ``changes of variable'' preserve the general form of $P$? Here is an answer:

\begin{lemma}\label{chvar} For $g\in K^\times$ and $\phi:= g^{-2}$ we have
$$ g^3P_{\times g}^\phi(Y)\ =\ 4Y'' + g^3P(g)Y.$$
\end{lemma} 
\begin{proof} Let $g,\phi\in K^\times$. Then 
\begin{align*} P_{\times g}(Y)\ &=\ 4gY'' + 8g'Y' + (4g'' + fg)Y\ =\ 4gY'' + 8g'Y' + P(g)Y, \quad \text{so}\\
 P_{\times g}^\phi(Y)\ &=\ 4g(\phi^2 Y'' + \phi'Y') + 8g'\phi Y' + P(g)Y\\ 
 &=\ 4g\phi^2 Y'' + (4g\phi' +8g'\phi)Y' + P(g)Y.\end{align*}
Now $4g\phi' + 8g'\phi=0$ is equivalent to $\phi^{\dagger}=-2g^\dagger$, which holds for $\phi=g^{-2}$.
For this~$\phi$ we get $P_{\times g}^\phi(Y)=
g^{-3}\big(4Y'' + g^3P(g)Y\big)$, that is, $g^3P_{\times g}^\phi(Y) = 4Y'' + g^3P(g)Y$. 
\end{proof}

\subsection*{Hardy fields} A {\em Hardy field\/} is a subfield of 
$\Gi$ that is closed under the derivation of $\Gi$. A Hardy field $H$ is considered as an ordered valued differential field in the obvious way, and has 
$\R\cap H$ as its field of constants. Hardy fields are pre-$H$-fields, and $H$-fields if they contain $\R$. Here are some well-known extension results:

\begin{prop}\label{prop:Hardy field exts} Any Hardy field $H$ has the following Hardy field extensions: \begin{enumerate}
\item[(i)] $H(\R)$, the subfield of $\Gi$ generated by $H$ and $\R$;
\item[(ii)] $H^{\operatorname{rc}}$, the real closure of $H$ as defined in Lemma~\ref{b1};
\item[(iii)] $H(\ex^f)$ for any $f\in H$;
\item[(iv)] $H(f)$ for any $f\in \Go$ with $f'\in H$;
\item[(v)] $H(\log f)$ for any $f\in H^{>}$.
\end{enumerate}
If $H$ is contained in $\Ginf$, then
so are the Hardy fields in {\rm (i), (ii), (iii), (iv), (v)}; likewise with $\Gom$ instead of $\Ginf$. 
\end{prop}

\noindent
Note that (v) is a special case of (iv), since $(\log f)'=f'/f\in H$
for $f\in H^{>}$. Another special case of (iv) is that $H(x)$ is a Hardy field. A consequence of the Proposition is that any Hardy field $H$ has a smallest real closed Hardy field extension $H^*$ with~$\R\subseteq H^*$ such that for all $f\in H^*$ we have $\ex^f\in H^*$
and $g'=f$ for some~$g\in H^*$. Note that then $H^*$ is Liouville
closed as defined in \cite[Section~10.6]{ADH}. 

We also have the following more general extension result from Rosenlicht~\cite{Ros}, attributed there to M. Singer: 

\begin{prop}\label{singer} Let $H$ be a Hardy field and $p(Y),q(Y)\in H[Y]$. Suppose $f\in \mathcal{C}^1$ is a solution of the differential equation
$y'q(y)=p(y)$ and $q(f)$ is a unit of $\mathcal{C}^1$. Then~$f$ generates a Hardy field $H(f)$ over $H$.
\end{prop}

\subsection*{Compositional inversion and compositional conjugation in Hardy fields} Let $H$ be a Hardy field, and
let $g\in \C^1$ be such that $g>\R$
and $g'\in H$. Then we have a Hardy field $H(g)$, and the compositional inverse 
$g^{\inv}\in \C^1$ of $g$ satisfies 
$$g^{\inv}>\R, \qquad (g^{\inv})'\ =\ (1/g')\circ g^{\inv}\in H\circ g^{\inv}$$
and yields an ordered field isomorphism
$$h\mapsto h\circ g^{\inv}\ :\ H \to H\circ g^{\inv}$$
such that for all $h\in H$,
$$(h\circ g^{\inv})'\ =\ (h'\circ g^{\inv})\cdot (g^{\inv})'\ =\ (h'/g')\circ g^{\inv}\in H\circ g^{\inv}.$$
Thus $H\circ g^{\inv}$ is again a Hardy field, and for $\phi=g'$
this yields an isomorphism
$$  h\mapsto h\circ g^{\inv}\ :\ H^{\phi} \to H\circ g^{\inv}$$
of pre-$H$-fields. 
If $H\subseteq \Ginf$ and $g\in \Ginf$, then
$H\circ g^{\inv}\subseteq \Ginf$; likewise with $\Gom$ instead of $\Ginf$.  For later use, a {\em $\C^{\infty}$-Hardy field\/} is a Hardy field $H\subseteq \Ginf$, and a {\em $\C^{\omega}$-Hardy field\/} (also called an {\em analytic Hardy field\/}) is a Hardy field $H\subseteq \Gom$.


\section{Extending Hardy Fields to $\upo$-free Hardy Fields}\label{secupo} 

\noindent
In this section we assume familiarity with \cite[Sections~5.2, 11.5--11.8]{ADH}. Here we summarize some of this material, and then use this to prove Theorem~\ref{upo} below. In the {\em Notations and terminology} at the end of the introduction we defined for any differential ring $R$ functions $\omega\colon R\to R$ and $\sigma\colon R^\times\to R$.  We define likewise 
$$\omega\ :\ \C^1_a\to \C^0_a, \qquad \sigma\ :\ (\C^1_a)^\times \to \C^0_a$$ 
by
$$\omega(z)\ =\ -2z'-z^2\quad\text{ and }\quad\sigma(y)\ =\ \omega(z)+y^2\text{ for $z:= -y^\dagger$.}$$
To clarify this
role of $\omega$ and $\sigma$ in connection with second-order linear differential equations, let $f\in \C^0_a$ and consider the differential equation 
$$4Y''+fY\ =\ 0.$$ Suppose $y\in \C^2_a$ is a {\em non-oscillating} solution, that is, a solution with $y(t)\ne 0$
for all sufficiently large $t$, say for all $t\ge b$, where 
$b\ge a$.
Then $z\in \C^1_b$ given by
$z(t)=2y'(t)/y(t)$ satisfies the first-order differential equation
$-2z'-z^2=f$ on $[b,\infty)$. Thus if the germ of $f$ at $+\infty$ belongs to a Hardy field $H$, then by Proposition~\ref{singer} the germ of $z$ at $+\infty$ (also denoted by $z$) generates a Hardy field~$H(z)$
with $\omega(z)=f$, which in turn
yields a Hardy field $H(z,y)$ with $y$ now denoting its germ at $+\infty$ so that $y\in (\Gi)^\times$ and $2y^\dagger=z$ in $\Gi$.
Thus $y_1:=y$ lies in a Hardy field extension of $H$. 
The germ $y_2$ of the function $[b,+\infty)\to\R$ given by $t\mapsto y(t)\int_b^t \frac{1}{y(s)^2}\,ds$
also satisfies $4y_2''+fy_2=0$, and $y_1$, $y_2$ are $\R$-linearly independent~\cite[Chapter~6, Lemma~3]{Bellman}. By Proposition~\ref{prop:Hardy field exts}(iv), $y_2$ lies in a Hardy field
extension of $H\<y_1\>=H(y,z)$; see also \cite[Theorem~2, Corollary~2]{Ros}. 

There might not exist a non-oscillating solution $y$, but we do have $\R$-linearly independent solutions $y_1, y_2\in \C^2_a$.
We saw before that 
$w:= y_1y'_2-y'_1y_2\in \R^\times$. 
Set $y:= y_1+ y_2\imag$. Then $4y'' + fy=0$ and 
$y(t)\ne 0$ for all $t\ge a$, and for $z\in \C^1_a[\imag]$ given by
$z(t)=2y'(t)/y(t)$ we have $-2z'-z^2=f$. Now
\begin{align*} z\ =\ \frac{2y_1'+2\imag y_2'}{y_1+\imag y_2}\ &=\
\frac{2y_1'y_1+2y_2'y_2- 2\imag(y_1'y_2-y_1y_2')}{y_1^2+y_2^2}\ =\ \frac{2(y_1'y_1+y_2'y_2)+2\imag w}{y_1^2+y_2^2},\\
\text{so }\ \operatorname{Re}(z)\ &=\ \frac{2(y_1'y_1+y_2'y_2)}{y_1^2+y_2^2}\in \C^1_a, \qquad \operatorname{Im}(z)\ =\ \frac{2w}{y_1^2+y_2^2}\in \C^2_a.
\end{align*} 
Thus $\operatorname{Im}(z)\in (\C^2_a)^\times$ and $\operatorname{Im}(z)^\dagger=-\operatorname{Re}(z)$  and 
$\sigma\big(\operatorname{Im}(z)\big)=
\omega(z)=f$ in $\C^1_a$.
Replacing $y_1$ by $-y_1$ changes 
$w$ to $-w$; in this way we can arrange that~${w>0}$.

\subsection*{The property of $\upo$-freeness}
Let $H\supseteq \R$ be a Liouville closed Hardy field.
Note that then $x\in H$ and $\log f\in H$ for all $f\in H^{>}$.  To express the property of $\upo$-freeness for $H$ we introduce the ``iterated logarithms'' $\ell_{\rho}$; more precisely, transfinite recursion yields a  sequence
$(\ell_\rho)$ in $H^{>\R}$ indexed by the ordinals $\rho$ less than some infinite limit ordinal $\kappa$ as follows:  $\ell_0=x$, and
$\ell_{\rho+1}:=\log \ell_\rho$; if $\lambda$ is an infinite limit ordinal
such that all $\ell_\rho$ with $\rho<\lambda$ have already been chosen,
then we pick $\ell_\lambda$ to be any element in $H^{>\R}$ such that $\ell_\lambda\prec \ell_\rho$ for
all $\rho<\lambda$, if there is such an $\ell_\lambda$, while if there is no
such $\ell_\lambda$, we put $\kappa:=\lambda$. From $(\ell_\rho)$ we obtain the sequences~$(\upg_\rho)$ in~$H^{>}$ and~$(\upl_\rho)$ in $H$ as follows:
$$
\upg_\rho\ :=\ \ell_{\rho}^\dagger, \qquad
\upl_\rho\ :=\  -\upg_\rho^\dagger\ =\ -\ell_\rho^\dagger{}^\dagger\ :=\ -(\ell_\rho^\dagger{}^\dagger).$$
Then $\upl_{\rho+1}=\upl_{\rho}+\upg_{\rho+1}$ and we have
\begin{align*} \upg_0\ &=\ \ell_0^{-1}, & \upg_1\ &=\ (\ell_0\ell_1)^{-1}, & \upg_2\ &=\ (\ell_0\ell_1\ell_2)^{-1},\\
\upl_0\ &=\ \ell_0^{-1}, & \upl_1\ &=\ \ell_0^{-1} + (\ell_0\ell_1)^{-1}, &
  \upl_2\ &=\ \ell_0^{-1} + (\ell_0\ell_1)^{-1} + (\ell_0\ell_1\ell_2)^{-1},
\end{align*}
and so on. Indeed, $v(\upg_\rho)$ is strictly increasing as a function of $\rho$
and is cofinal in~$\Psi_H=\big\{v(f^\dagger):f\in H,\ 0\neq f\nasymp 1\big\}$; we refer to
\cite[Section~11.5]{ADH} for this and some of what follows.
Also, $(\upl_\rho)$ is a strictly increasing pc-sequence which is cofinal in $\Upl(H)$; see \cite[Section~11.8]{ADH} for the definition of the set $\Upl(H)\subseteq H$, which is downward closed since $H$ is Liouville closed. 
The latter also gives that $H$ is
$\upl$-free as defined in \cite[Section~11.6]{ADH}, equivalently, 
$(\upl_{\rho})$ has no pseudolimit in $H$. The function $\omega\colon H \to H$ is strictly increasing on~$\Upl(H)$ and setting $\upo_{\rho}:= \omega(\upl_{\rho})$ we obtain a strictly increasing
pc-sequence $(\upo_{\rho})$ which is cofinal in 
$\omega\big(\Upl(H)\big)=\omega(H)$:
$$ \upo_0\ =\ \ell_0^{-2},\qquad \upo_1\ =\ \ell_0^{-2} + (\ell_0\ell_1)^{-2},
\qquad \upo_2\ =\ \ell_0^{-2} + (\ell_0\ell_1)^{-2} + (\ell_0\ell_1\ell_2)^{-2},$$
and so on; see \cite[Sections~11.7,~11.8]{ADH} for this and some of what follows. Now $H$ being $\upo$-free is equivalent to $(\upo_\rho)$ having no pseudolimit in $H$. By 
\cite[Corollary~11.8.30]{ADH} the pseudolimits of $(\upo_{\rho})$
in $H$ are exactly the $\upo\in H$ such that $\omega(H) < \upo < \sigma\big(\Upg(H)\big)$. Here the upward closed subset $\Upg(H)$ of $H$ is given by
$$\Upg(H)\ =\ \big\{a^\dagger:\ a\in H,\ a\succ 1\big\}\ =\ \{a\in H:\ a>\upg_{\rho} \text{ for some $\rho$}\},$$
and $\sigma$ is strictly increasing on $\Upg(H)$. Thus $H$ is not
$\upo$-free if and only if there exists an $\upo\in H$ such that $\omega(H) < \upo < \sigma\big(\Upg(H)\big)$.

\medskip\noindent
We are now ready to prove the following:

\begin{theorem}\label{upo} Every Hardy field has an $\upo$-free Hardy field extension.
\end{theorem} 
\begin{proof} It is enough to show that every maximal Hardy field is $\upo$-free. That reduces to showing that
every non-$\upo$-free Liouville closed Hardy field containing $\R$
has a proper Hardy field extension. So assume $H\supseteq \R$ is a Liouville closed Hardy field and $H$ is not $\upo$-free. We shall construct a proper Hardy field extension
of $H$. We have $\upo\in H$ such that 
$$\omega(H)\ <\  \upo\  <\ \sigma\big(\Upg(H)\big).$$ 
Take $a\in \R$ such that $\upo$ is the germ of a function in 
$\C^2_a$, this function also to be denoted by $\upo$.
With $\upo$ in the role of $f$ in the discussion preceding the statement of the theorem, we have  $\R$-linearly independent
solutions $y_1, y_2\in \C^2_a$ of the differential equation 
$4Y'' + \upo Y=0$ whose germs at $+\infty$ (also denoted by $y_1$ and $y_2$)
lie in~$\Gi$. Then the complex solution $y=y_1+y_2\imag$ is a unit
of $\C^2_a[\imag]$, and so we have 
$z:=2y^\dagger\in \C^1_a[\imag]$. The germs of $y$ and $z$ at $+\infty$ are also denoted by $y$ and $z$ and lie in $\Gi[\imag]$. We shall prove that the elements $\operatorname{Re}(z)$ and $\operatorname{Im}(z)$ of $\Gi$ generate a Hardy field extension $K=H\big(\operatorname{Re}(z), \operatorname{Im}(z)\big)$ of $H$
with $\upo=\sigma\big(\operatorname{Im}(z)\big)\in \sigma(K^\times)$. 
We can assume that $w:= y_1y_2'-y_1'y_2\in \R^{>}$, so $\operatorname{Im}(z)(t)>0$ for all $t\ge a$.

We have $\upo_{\rho} \leadsto \upo$, with $\upo-\upo_{\rho}\sim \upg_{\rho+1}^2$ by \cite[Lemma~11.7.1]{ADH}.
We set
$g_{\rho}:=\upg_{\rho}^{-1/2}$, so $2g_{\rho}^\dagger=\upl_{\rho}=-\upg_{\rho}^\dagger$. 
For $h\in H^\times$ we also have $\omega(2h^\dagger)=-4h''/h$,
hence $P:= 4Y'' + \upo Y\in H\{Y\}$ gives
$$P(g_{\rho})\ =\  
g_{\rho}(\upo-\upo_{\rho})\ \sim\ g_{\rho}\upg_{\rho+1}^2,$$ and so with an eye towards using Lemma~\ref{chvar}:
$$g_{\rho}^3P(g_{\rho})\ \sim\ g_{\rho}^4\upg_{\rho+1}^2\ \sim\ \upg_{\rho+1}^2/\upg_{\rho}^2\ \asymp\ 1/\ell_{\rho+1}^2.$$
Thus with
$g:= g_{\rho}=\upg_{\rho}^{-1/2}$, $\phi=g^{-2}=\upg_{\rho}$ we  have $A_{\rho}\in \R^{>}$ such that
\begin{equation}\label{eq:bound}
g^3P_{\times g}^{\phi}(Y)\ =\ 4Y'' + g^3P(g)Y,\quad 
|g^3P(g)(t)|\ \le\  A_{\rho}/\ell_{\rho+1}(t)^2, \text{ eventually.}
\end{equation}
From  $P(y)=0$ we get $P^{\phi}_{\times g}(y/g)=0$, that is,
$y/g\in \Gi[\imag]^\phi$ is a solution of $4Y'' + g^3P(g)Y=0$, with $g^3P(g)\in H\subseteq \Gi$. 
Now $\ell_{\rho+1}'=\ell_{\rho}^\dagger=\phi$, so the end of the previous section yields the isomorphism
$H^{\phi} \to H\circ \ell_{\rho+1}^{\inv}$ of $H$-fields, where~$\ell_{\rho+1}^{\inv}$ is the compositional inverse of 
$\ell_{\rho+1}$. Under this isomorphism the equation $4Y''+g^3P(g)Y=0$ corresponds to the equation
$$4Y'' + f_{\rho}Y\ =\ 0, \qquad f_{\rho}\ :=\ g^3P(g)\circ \ell_{\rho+1}^{\inv}\in H\circ \ell_{\rho+1}^{\inv}\ \subseteq\ \Gi.$$ 
The equation $4Y'' + f_{\rho}Y=0$ has the ``real'' solutions
$$y_{i,\rho}\ :=\ (y_i/g)\circ \ell_{\rho+1}^{\inv}\in \Gi\circ \ell_{\rho+1}^{\inv}\ =\ \Gi \qquad(i=1,2),$$  and the ``complex''
solution $$y_\rho:=y_{1,\rho} + y_{2,\rho}\imag=(y/g)\circ \ell_{\rho+1}^{\inv},$$ which is a unit of the ring
$\Gi[\imag]$. We set $z_{\rho}:=2 y_{\rho}^\dagger\in \Gi[\imag]$. 
The bound in \eqref{eq:bound} gives 
$$|f_{\rho}(t)|\ \le\  A_{\rho}/t^2, \text{ eventually},$$
which by Corollary~\ref{cor:bound} yields  positive constants $B_{\rho}$, $c_\rho$ such that
$$|z_{\rho}(t)|\ \le\ B_{\rho}t^{c_\rho}, \quad \text{ eventually.}$$ Using $(\ell_{\rho+1}^{\inv})'=(1/\ell_{\rho+1}')\circ \ell_{\rho+1}^{\inv}$ we obtain
$$z_{\rho}\ =\ 2\big((y/g)^\dagger\circ \ell_{\rho+1}^{\inv}\big)\cdot (\ell_{\rho+1}^{\inv})'\ =\ 2\big((y/g)^\dagger/\ell_{\rho+1}'\big)\circ \ell_{\rho +1}^{\inv}\ =\ \big((z-2g^\dagger)/\ell_{\rho+1}'\big)\circ \ell_{\rho+1}^{\inv}.$$ 
In combination with the eventual bound on $|z_\rho(t)|$ this yields
\begin{align*} \left|\frac{z(t)-2g^\dagger(t)}{\ell_{\rho+1}'(t)}\right|\ &\le\ B_{\rho}\,\ell_{\rho+1}(t)^{c_\rho}\ \text{ eventually, hence}\\
  |z(t)- \upl_{\rho}(t)|\ &\le\ B_{\rho}\,\ell_{\rho+1}(t)^{c_\rho}\,\ell_{\rho+1}'(t)\ =\ 
  B_{\rho}\,\ell_{\rho+1}(t)^{c_\rho}\,\upg_{\rho}(t), \text{ eventually, so}\\ 
z(t)\ &=\ \upl_{\rho}(t)  + R_{\rho}(t), \quad |R_{\rho}(t)|\le B_{\rho}\,\ell_{\rho+1}(t)^{c_\rho} \,\upg_{\rho}(t), \text{ eventually}.
\end{align*}
We now use this last estimate with $\rho+1$ instead of $\rho$, 
together with 
$$\upl_{\rho+1}\ =\ \upl_{\rho}+\upg_{\rho+1},\quad  \ell_{\rho+1}\upg_{\rho+1}\ =\ \upg_{\rho}.$$
This yields
\begin{align*} z(t)\, &=\, \upl_{\rho}(t) + \upg_{\rho+1}(t) + R_{\rho+1}(t)\ \text{ eventually, with}\\
 |R_{\rho+1}(t)|\, &\le\, B_{\rho+1}\,\ell_{\rho+2}(t)^{c_{\rho+1}}\, \upg_{\rho+1}(t) \\ &\, =\, B_{\rho+1}\big(\ell_{\rho+2}(t)^{c_{\rho+1}}/\ell_{\rho+1}(t)\big)\,\upg_{\rho}(t)\, \text{ eventually},\\
\text{so }z(t)\, &=\, \upl_{\rho}(t) + o\big(\upg_{\rho}(t)\big)\ \text{ as }\ t\to \infty,\ \text{ and thus}\\
\operatorname{Re}(z)(t)\, &=\, \upl_{\rho}(t) + o\big(\upg_{\rho}(t)\big), \quad\operatorname{Im}(z)(t)\, =\, o\big(\upg_{\rho}(t)\big), \text{ as }\ t\to \infty.
\end{align*}
Recall that $(\upl_{\rho})$ is a strictly increasing divergent
pc-sequence $(\upl_{\rho})$ in $H$ which is cofinal in $\Upl(H)$. By the above,
$\upl:= \operatorname{Re}(z)\in \Gi$ satisfies $\Upl(H) < \upl < \Upd(H)$. This yields an ordered subfield 
$H(\upl)$ of $\Gi$, which by Lemma~\ref{ps1} is an immediate
valued field extension of $H$ with $\upl_{\rho} \leadsto \upl$.

Pick functions in $\C^0_a$ whose germs at $+\infty$ are the elements
$\ell_{\rho}$,~$\upg_{\rho}$,~$\upl_{\rho}$ of $H$; we denote these functions also by $\ell_{\rho}$,~$\upl_{\rho}$,~$\upg_{\rho}$. From 
$\ell_{\rho}^\dagger=\upg_{\rho}$ and $\upg_{\rho}^\dagger=-\upl_{\rho}$ in $H$ we obtain constants $c_{\rho}, d_{\rho}\in \R^{>}$ such that for all $t\ge a$,
$$ \ell_{\rho}(t)\ =\ c_{\rho}\exp\left[\int_a^t \upg(s)\,ds\right], \quad \upg_{\rho}(t)\ =\ d_{\rho}\exp\left[-\int_a^t\upl_{\rho}(s)\,ds\right]. $$
Set $\upg:= \operatorname{Im}(z)$, so $\upg^\dagger=-\upl$, and both
$\upg$ and $\upl$ are already given as elements of $\C^0_a$.
Since $\upg(t)>0$ for all $t\ge a$ we have a constant $d\in \R^{>}$ such that for all $t\ge a$, 
$$ \upg(t)\ =\ d\exp\left[-\int_a^t\upl(s)\,ds\right].$$ 
The above estimate for $\upl=\operatorname{Re}(z)$ gives 
$$\upl_{\rho}(t)\ <\ \upl(t)\ <\ \upl_{\rho}(t) + \upg_{\rho}(t), \text{ eventually,}$$
so we have constants $a_{\rho}, b_{\rho}\in \R$ such that
$$\int_a^t \upl_{\rho}(s)\,ds\ <\ a_{\rho}+ \int_a^t \upl(s)\,ds\ <\ b_{\rho} +\int_a^t \upl_{\rho}(s)\,ds + \int_a^t \upg_{\rho}(s)\,ds,\ \text{ eventually},$$
which by applying $\exp(-*)$ yields 
$$\frac{1}{d_{\rho}}\upg_{\rho}(t)\ >\ \frac{1}{\ex^{a_{\rho}}d}\upg(t)\ >\ \frac{c_{\rho}}{\ex^{b_{\rho}}d_{\rho}} \upg_{\rho}(t)/\ell_{\rho}(t),\
 \text{ eventually}.$$
Here the positive constant factors don't matter,
since the valuation of $\upg_{\rho}$ is strictly increasing and that of $\upg_{\rho}/\ell_{\rho}=(1/\ell_{\rho})'$
is strictly decreasing with $\rho$. Thus 
for all~$\rho$ we have 
$\upg_{\rho} > \upg > (1/\ell_{\rho})'$, in $\Gi$. In view of Lemma~\ref{ps2} applied to $H(\upl)$,~$\upg$ in the role of
$K$,~$f$ this yields an ordered subfield $H(\upl, \upg)$
of $\Gi$. Moreover, $\upg$ is transcendental over $H(\upl)$ with
$\upg^\dagger= -\upl$, and $\upg$ 
satisfies the second-order differential equation
$2yy''-3(y')^2+y^4-\upo y^2=0$ over $H$ (obtained from the relation 
$\sigma(\upg)=\upo$ by multiplication with $\upg^2$). It follows that $H(\upl,\upg)$ is closed under the derivation of~$\Gi$, and hence $H(\upl, \upg)=H\<\upl\>$ is a Hardy field. 
\end{proof} 

\noindent
The proof also shows that every $\C^\infty$-Hardy field
has an 
$\upo$-free $\C^\infty$-Hardy field extension, and the same with 
$\C^\omega$ instead of $\C^\infty$.

\end{document}